\documentclass[11pt,centertags]{amsart}
\usepackage{amssymb}            
\usepackage{mathrsfs}             
\usepackage{hyperref}           
\usepackage{ifsym}              
\usepackage{calc}               
\usepackage{enumerate}          
\usepackage{tikz-cd}			
\usepackage[all]{xy}
\usepackage{psfrag}             
\usepackage{verbatim}           
\usepackage{fancyhdr}           
\usepackage{url}                
\usepackage{gloss}              
\usepackage{yhmath}             

\usepackage{ifpdf}

\newtheorem*{main*}{Main Theorem}

\newtheorem{theorem}{Theorem}[section]
\newtheorem*{theorem*}{Theorem}
\newtheorem{proposition}[theorem]{Proposition}
\newtheorem{lemma}[theorem]{Lemma}

\newtheorem*{question*}{Question}

\newtheorem*{conjecture*}{Conjecture}

\theoremstyle{definition}

\newtheorem{definition}[theorem]{Definition}
\newtheorem*{definition*}{Definition}

\newtheorem*{example*}{Example}

\theoremstyle{remark}
\newtheorem{remark}[theorem]{Remark}

\numberwithin{equation}{section}

\DeclareMathOperator{\Sp}{\Sp}

		\renewcommand{\bf}{\bfseries}

		\renewcommand{\it}{\itshape}

		\renewcommand{\hat}{\widehat}

\begin{document}

\author[Jean-Fran\c{c}ois Lafont]{Jean-Fran\c{c}ois Lafont$^\dagger$}
\author[Gangotryi Sorcar]{Gangotryi Sorcar}
\author[Fangyang Zheng]{Fangyang Zheng$^*$}

\thanks{$^\dagger$ J.-F.L. is partly supported by the NSF, under grant DMS-1510640.}
\thanks{$^*$ F.Z. is partially supported by a Simons Collaboration Grant from the Simons Foundation.}

\title{Some K\"ahler structures on products of $2$-spheres}

\address{Ohio State University, Columbus, OH 43210, USA} 
\email{lafont.1@osu.edu}

\address{Ohio State University, Columbus, OH 43210, USA} 
\email{sorcar.1@osu.edu} 

\address{Ohio State University, Columbus, OH 43210, USA} 
\email{zheng.31@osu.edu}



\begin{abstract}
We consider a family of K\"ahler structures on products of $2$-spheres, arising from complex Bott manifolds.
These are obtained via iterated $\mathbb P^1$-bundle constructions, generalizing the classical Hirzebruch
surfaces. We show that the resulting K\"ahler structures all have identical Chern classes. We construct Bott
diagrams, which are rooted forests with an edge labelling by positive integers, and show that these classify 
these K\"ahler structures up to biholomorphism. 
\end{abstract}

\maketitle


\thispagestyle{empty} 

\section{Introduction}

In complex geometry, it is interesting to study the class of complex structures (or K\"ahler structures) supported on a fixed smooth
oriented manifold $M$. Since the basic invariants of a complex manifolds are the Chern classes, it is tempting to try and use these
to distinguish complex structures on $M$. In complex dimension two, the Hirzebruch surfaces $\mathbb F _m$
are topologically either diffeomorphic to $S^2 \times S^2$ (if $m$ is even), or to $\mathbb P^2 \# \overline{\mathbb P^2}$ (if $m$ is odd).
Focusing on the Hirzebruch surfaces diffeomorphic to $S^2 \times S^2$, a celebrated result of Hirzebruch \cite{hirzebruch}
shows that all the $\mathbb F_{2k}$ 
are distinct as complex manifolds, even though they have identical Chern classes.

In the present paper, we extend Hirzebruch's result, by considering {\it $\mathbb Z$-trivial complex Bott manifolds} (see Bott and
Samelson \cite{BS}). In complex dimension
two, these are precisely the Hirzebruch surfaces $\mathbb F _{2k}$. In complex dimension $n$, these are compact K\"ahler manifolds 
diffeomorphic to $S^2 \times \cdots \times S^2 = (S^2)^n$.
To each $n$-dimensional $\mathbb Z$-trivial complex Bott manifold $M$, 
we associate a {\it Bott diagram}, which is a rooted forest equipped with an edge labelling by positive integers. Our main result is the following:

\begin{main*}
Every $n$-vertex rooted forest equipped with an edge labelling by positive integers arises as the Bott diagram of some
$n$-dimensional $\mathbb Z$-trivial complex Bott manifold. Moreover, for an arbitrary pair of $n$-dimensional $\mathbb Z$-trivial 
complex Bott manifolds $M_1, M_2$, we have:
\begin{enumerate}
\item $M_1$ is biholomorphic to $M_2$ if and only if their Bott diagrams are isomorphic, but
\item there is a diffeomorphism $\phi: M_1\rightarrow M_2$ with the property that $\phi^*(c(M_2))=c(M_1)$, where $c$ 
denotes the total Chern classes.
\end{enumerate}
\end{main*}

Our result provides a combinatorial classification of a certain family of K\"ahler structures on the products 
$S^2 \times \cdots \times S^2 = (S^2)^n$. When $n=2$, the only K\"ahler structures on $S^2\times S^2$ 
are those arising from the Hirzebruch surfaces. When $n\geq 3$, we do not know whether these products of 
$2$-spheres support any other K\"ahler structures. Our result also shows that these K\"ahler structures are
indistinguishable as far as Chern classes are concerned.

\vskip 5pt

Recall that the rational Pontrjagin classes of a smooth manifold are defined using the smooth structure. A celebrated result
of Novikov \cite{Novikov} shows however that these classes in fact only depend on the underlying topological structure (other
proofs were given in \cite{Gromov}, \cite{Ranicki}, \cite{RY}, and \cite{RW}). More precisely, if one 
has a pair of homeomorphic smooth manifolds, then the homeomorphism can be chosen to take the total rational Pontrjagin class 
to the total rational Pontrjagin class. Because of the lack of known counterexamples, we would like to raise the following analogous 
question for Chern classes:

\begin{question*}
Let $X$ be a compact smooth manifold and let ${\mathcal M}_0$ be the set of biholomorphism classes of all K\"ahlerian complex structures 
on $X$.  Do the total Chern classes $c(J)\in H^{2\ast }(M, {\mathbb Z})$ (for all $J\in {\mathcal M}_0$) belong to a single orbit for the action
of the diffeomorphism group on the cohomology ring? 
\end{question*}

Part of the difficulty in addressing this question is the lack of examples.  Indeed, several classes of manifolds are known to support 
unique K\"ahler structures -- see for instance the rigidity results of Hirzebruch-Kodaira \cite{HK}, Yau \cite{yau} and Siu \cite{siu}. 
In contrast, there are very few smooth manifolds that are known to support multiple distinct K\"ahler structures. 
Our main theorem shows that the answer to this question is {\it yes} for a particular family of K\"ahler structures
on $X= (S^2)^n$.  

We also point out that there are examples of manifolds where the answer to this question is {\it no}. Notably, it follows 
from work of Kotschick \cite{Kotschick1} \cite{Kotschick2}, that there exist a pair of K\"ahler $3$-folds $M_1, M_2$ which are diffeomorphic,
but have distinct Chern numbers $c_1^3(M_1)\neq c_1^3(M_2)$. A consequence is that their first Chern classes must lie in distinct orbits of
the diffeomorphism group action on second cohomology. Similar examples exist in higher dimensions. Nevertheless, the question might
have an affirmative answer for some manifolds, particularly in the presence of a high degree of symmetry.

\vskip 5pt

Our paper is structured as follows. We review some background material in Section \ref{sec-backgrnd}, and prove our 
main theorem in Section \ref{sec-main-thm}. Our argument requires a reconstruction result for labelled rooted forests, 
which is explained in Section \ref{sec-reconstruction}.


\section{Background material}\label{sec-backgrnd}


\subsection*{Bott manifolds.}
Recall that a Bott manifold $M^n$ is one that admits a {\em Bott Tower,} namely, $M^n=B_n$ and
\begin{equation}
\label{tower} B_n \xrightarrow {\pi_n} B_{n\!-\!1} \xrightarrow {\pi_{n\!-\!1}} \cdots \xrightarrow {\pi_2} B_1 \xrightarrow{\pi_1} B_0=\{ \text{a  point} \}
\end{equation}
where for each $1\leq j\leq n$, $B_j={\mathbb P}({\mathcal O}\oplus S_j)$ is the projectivization of the direct sum of the trivial line bundle with a holomorphic line bundle $S_j$ over $B_{j\!-\!1}$, with $\pi_j$ the projection map. 

Clearly, $B_1={\mathbb P}^1$, and $B_2$ is a Hirzebruch surface ${\mathbb F}_m = {\mathbb P}({\mathcal O}_{{\mathbb P}^1} \oplus {\mathcal O}_{{\mathbb P}^1}(-m))$ over ${\mathbb P}^1$, where $m$ is any nonnegative integer. It is well-known that all ${\mathbb F}_{2k}$ are diffeomorphic to ${\mathbb F}_0=\mathbb P^1 \times \mathbb P^1 \cong S^2\!\times \!S^2$, while all ${\mathbb F}_{2k-1}$ are diffeomorphic to ${\mathbb F}_1={\mathbb P}^2 \# \overline{{\mathbb P}^2}$, the one point blow up of $ {\mathbb P}^2$.

\begin{definition}
A Bott manifold $M^n$ is called ${\mathbb Z}$-trivial, if it is diffeomorphic to $({\mathbb P}^1)^n$, the product of $n$ copies of the complex projective line.
\end{definition}

By the work of Choi and Masuda \cite{ChoiMasuda}, a Bott manifold $M^n$ is ${\mathbb Z}$-trivial if its integral cohomology ring is isomorphic (as a graded ring) to that of $({\mathbb P}^1)^n$. In fact, they show that every graded ring isomorphism between the cohomology rings of two ${\mathbb Q}$-trivial Bott manifolds is induced by a diffeomorphism. Here $M^n$ is ${\mathbb Q}$-trivial means that the cohomology ring of $M$ with ${\mathbb Q}$-coefficients is isomorphic (as a graded ring) to that of $({\mathbb P}^1)^n$.


\subsection*{Projectivization of vector bundles.} Let us recall some general facts concerning projectivizations of vector bundles. 

Let $E$ be a holomorphic vector bundle of rank $r$ over a compact complex manifold $B$, let $\pi : M={\mathbb P}(E)\rightarrow B$ 
be the  projectivization of $E$, where $\pi$ is the projection map. We adopt the algebro-geometric convention here, namely,  
$\pi^{-1}(x)={\mathbb P}(E_x)$ is the set of all the {\em hyperplanes} (instead of {\em lines}) through the origin in the fiber 
$E_x\cong {\mathbb C}^{r\!-\!1}$.  $M$ is again a compact complex manifold, a holomorphic fiber bundle with fiber 
${\mathbb P}^{r\!-\!1}$ over $B$. 

Denote by $L$ the dual of the tautological line bundle, then we have the following two short exact sequences of holomorphic vector bundles 
over $M$:
\begin{eqnarray}
\label{euler}  &&  0 \rightarrow  {\mathcal O}_M \rightarrow  \pi^{\ast } E^{\ast } \otimes L \rightarrow  T_{M|B} \rightarrow  0 \\
 \label{relativetangent} &&  0 \rightarrow  T_{M|B} \rightarrow  T_{M} \rightarrow  \pi^{\ast } T_{B} \rightarrow  0 
\end{eqnarray}
where $T_{M|B}$ is the relative tangent bundle, namely, the kernel of the differential of $\pi$. The first Chern class $\xi = c_1(L)$ 
satisfies the Grothendieck equation
\begin{equation*}
 f(\xi) : = \xi^r - \xi^{r\!-\!1} \cdot \pi^{\ast }c_1(E) + \xi^{r\!-\!2} \cdot \pi^{\ast }c_2(E) - \cdots + (-1)^r \pi^{\ast }c_r(E) = 0,
\end{equation*}
while the cohomology ring (or the Chow ring) of $M$ is generated by the pull back of that of $B$ and $\xi$:
\begin{equation}
H^{\ast }(M, {\mathbb Z})  : = \pi^{\ast }H^{\ast }(B, {\mathbb Z})\  [\xi ]  \ / \  ( f(\xi )) .
\end{equation}

Recall that a {\em section} of $\pi$ is a complex submanifold $Z\subseteq M$ such that $\pi|_Z: Z \rightarrow B$ is a biholomorphism. 
Equivalently, a section of $\pi$ is given by a holomorphic map $i : B \rightarrow M$ such that $\pi \circ i = id_B$. In this case the image 
$i(B)$ is the submanifold in $M$ isomorphic to $B$. Note that the sections of $\pi$ correspond to quotient line bundles of $E$. 

To see this, let $Q$ be a holomorphic line bundle on $B$ which is a quotient bundle of $E$. As we are using the hyperplane convention for 
projectivizations, so ${\mathbb P}(Q)\cong B$ is a submanifold of ${\mathbb P}(E)=M$, which gives a section of $\pi$. Conversely, given a 
section $i: B\rightarrow M$ of $\pi$, since the tautological line bundle $L^{\ast}$ is a subbundle of $\pi^{\ast}E^{\ast }$ on $M$, $Q=i^{\ast }L$ 
would be a quotient line bundle of $i^{\ast} \pi^{\ast }E = E$ on $B$.

Next, let us specialize to the situation when the vector bundle on $B$ is $E={\mathcal O}\oplus S$, the sum of the trivial line bundle with 
another line bundle $S$. Writing $s= - \pi^{\ast}c_1(S)$, the above short exact sequences (\ref{euler}) , (\ref{relativetangent}), along with the 
Grothendieck equation, gives us
\begin{equation}
c_1(T_{M|B})=2\xi + s, \ \ c(M) = (1+ 2\xi +s ) \cdot \pi^{\ast }c(B), \ \ \mbox{and} \ \xi^2 = - \xi \cdot s
\end{equation}
in the cohomology (or the Chow) ring $H^{\ast }(M,{\mathbb Z})$. 


\subsection*{Cohomology ring of Bott manifolds.}

Now let us apply these formula to the $j$-th stage $\pi_j: B_j\rightarrow B_{j\!-\!1}$, which is the projectivization of the splitting rank $2$ vector bundle ${\mathcal O}\oplus S_j$ on $B_{j\!-\!1}$, then we get the following:
\begin{eqnarray*}
&& H^{\ast }(B_j,{\mathbb Z}) = \pi_j^{\ast } H^{\ast }(B_{j\!-\!1},{\mathbb Z}) [\xi_j] / (\xi_j^2+\xi_js_j)  \\  
&& c(B_j) = (1+ 2\xi_j + s_j ) \cdot \pi_j^{\ast }c(B_{j\!-\!1}),
\end{eqnarray*}
where $-s_j$ and $\xi_j$ are the first Chern class of $\pi_j^{\ast }S_j$ and $L_j={\mathcal O}_{\pi_j}(1)$, the dual of the tautological line bundle on $B_j$.

Given a Bott manifold $M^n$ with Bott tower \eqref{tower}, let us write
\begin{eqnarray*}
x_j & = &  (\pi_{j\!+\!1}\circ \cdots \circ \pi_n)^{\ast } \xi_j \\
h_j & = &  (\pi_{j\!+\!1}\circ \cdots \circ \pi_n)^{\ast } s_j
\end{eqnarray*}
for each $1\leq j\leq n$. Note that $x_1$ is the first Chern class of the pull back to $M$ of ${\mathcal O}_{{\mathbb P}^1} (1)$  on $B_1$, and $h_1=0$. By an inductive argument, we obtain the following:
\begin{eqnarray}
&& 
H^{\ast }(M, {\mathbb Z}) = {\mathbb Z} [ x_1, \ldots , x_n] / ( x_1^2, x_2^2+x_2h_2, \ldots , x_n^2+x_nh_n ) \\
&& c(M) = (1+2x_1)(1+2x_2+h_2) \cdots (1+2x_n+h_n)
\end{eqnarray}
where $x_1, \ldots , x_n$ is a set of generators for $H^2(M,{\mathbb Z})\cong {\mathbb Z}^n$, and each $h_j$ satisfies
\begin{equation}
h_j = a_{j1}x_1 + a_{j2}x_2 + \cdots + a_{j,j\!-\!1}x_{j\!-\!1} 
\end{equation}
where all $a_{jk}$ are integers. 

\begin{example*}
In the special case where all the line bundles $S_j$ are trivial, we get the product $P=({\mathbb P}^1)^n$ of $n$-copies of the complex projective line ${\mathbb P}^1$. In this case, all $h_j=0$ and we will denote the corresponding $x_j$ by $y_j$. The above computations give us:
\begin{eqnarray} 
&& H^{\ast }(P,{\mathbb Z}) = {\mathbb Z}[y_1, \ldots , y_n] / (y_1^2, \ldots , y_n^2) \\
&&  c(P)=(1+2y_1)\cdots (1+2y_n) 
\end{eqnarray}
\end{example*}

\vskip 10pt


\section{Proof of the Main Theorem.}\label{sec-main-thm}

This entire section is devoted to the proof of the Main Theorem.


 \subsection{The structure of ${\mathbb Z}$-trivial Bott manifolds.} We start by analyzing how the ${\mathbb Z}$-triviality condition
 affects the cohomology elements $h_j$.
 
For a given a Bott tower on $M^n$, assume that $2\mid h_j$ and $h_j^2=0$ for all $j$. Write $z_j=x_j + \frac{1}{2}h_j$. Since $h_1=0$, and 
for each $2\leq j\leq n$ the corresponding $h_j$ is generated by $x_1, \ldots , x_{j-1}$, it follows that $\{ z_1, \ldots , z_n\} $ generates 
$H^{\ast}(M,{\mathbb Z})$. Also, each $z_j^2 =0$ by the Grothendieck equation. So defining $\phi (y_j)=z_j$ gives a graded ring isomorphism 
$\phi : H^*(P; \mathbb Z) \rightarrow H^*(M^n ; \mathbb Z)$. By the result of Choi and Masuda \cite{ChoiMasuda}, there is a diffeomorphism
$\Phi: M^n\rightarrow P\cong (S^2)^n$, which induces $\Phi^* = \phi$. It follows that 
$M^n$ is ${\mathbb Z}$-trivial. Moreover, by the Chern class formula, we see that $\phi (c(P)) = c(M)$. 

Conversely, if there exists an isomorphism $\phi : H^{\ast}(P,{\mathbb Z}) \rightarrow H^{\ast}(M,{\mathbb Z})$, then we claim that 
$2\mid h_j$ and $h_j^2=0$ for all $j$. To see this, let us write $\phi (y_j)=z_j$. We have
$$ H^{\ast }(M,{\mathbb Z}) = {\mathbb Z} [z_1, \ldots , z_n] / ( z_1^2, \ldots , z_n^2).$$
For each $1\leq k\leq n$, the group $H^{2k}(M,{\mathbb Z})$ is a free abelian group generated by products
$z_I = z_{i_1}  \cdots z_{i_k}$ for all multi-indices $I= (i_1\cdots i_k)$ of length $k$, where $1\leq i_1 < i_2 < \cdots < i_k\leq n$. 
Note that for any integer linear combination $z=a_1z_1+ \cdots + a_nz_n$, if $z^2=0$, then $a_ia_j=0$ for all $i\neq j$, 
thus at most one of these $a_i$ could be non-zero.

Now we proceed to show that $2\mid h_j$ and $h_j^2=0$, by induction on $j$, where $j\in A:=\{ 1, 2, \ldots , n\}$. First we have $h_1=0$. 
For $j=2$, since $x_1^2=0$, we know that there must be a unique $i_1 \in A$ such that $x_1=\varepsilon_1z_{i_1}$, where 
$\varepsilon_1=\pm 1$ since $x_1$ is a primitive element in $H^2(M,{\mathbb Z})$. Write $x_2=az_{i_1}+z$, where $z$ is a 
linear combination of $z_j$ for $j\in A\setminus \{ i_1\}$. We have $h_2=bz_{i_1}$ since $h_2$ is a multiple of $x_1$. Since 
$x_2(x_2+h_2)=0$, we have
$$ (2a+b)z_{i_1}z + z^2 = 0. $$ 
Since $H^4(M, {\mathbb Z})$ is a free abelian group with generators $z_iz_j$ for $1\leq i<j\leq n$, we conclude from the above 
equality that $2a+b=0$ and $z^2=0$. So $2\mid h_2$, $h_2^2=0$, and $z=x_2+\frac{1}{2}h_2$ satisfies $z^2=0$, thus equals to 
$\varepsilon_2z_{i_2}$ for some $i_2\neq i_1$, and $\varepsilon_2=\pm 1$. 

Now assume that for a fixed $2\leq k<n$, we already have $2\mid h_j$, $h_j^2=0$ for each $j\leq k$, and 
$x_j':= x_j+\frac{1}{2}h_j=\varepsilon_j z_{i_j}$ where $i_1 $, ..., $i_k$ are all distinct in $A$ and $\varepsilon_j=\pm 1$. 
Since $h_{k+1}$ is a linear combination of $x_1', \ldots , x_k'$, we can write
$$ h_{k+1} = b_1z_{i_1} + \cdots + b_k z_{i_k} $$
Also, let us write $x_{k+1}=a_1z_{i_1} + \cdots + a_k z_{i_k} + z $, where $z$ is a combination of those $z_j$ for $j$ in 
$A\setminus \{ i_1, \ldots , i_k\}$. Now by applying the Grothendieck equation, namely, $x_{k+1}(x_{k+1}+ h_{k+1})=0$, we 
get the equation
$$ \sum_{j=1}^k (2a_j + b_j)z_{i_j}z + z^2 + \sum_{j,l=1}^k a_l(a_j+ b_j)z_{i_j} z_{i_l} =0. $$
Since $z$ cannot be zero, we know that $b_j=-2a_j$ for each $j\leq k$, so $2\mid h_{k+1}$ and $h_{k+1}^2=0$. Furthermore, 
$x_{k+1}'=z$ is a square zero primitive element, thus must be of the form $\pm z_{i_{k+1}}$ for some $i_{k+1}$ in 
$A\setminus \{ i_1, \ldots , i_k\}$.

To summarize, we have established the following (also independently obtained by J. H. Kim \cite{kim}):

\begin{lemma}\label{diffeo-orbit}
If $M^n$ is a Bott manifold and $\phi$ is an isomorphism between the integral cohomology rings of $P=({\mathbb P}^1)^n$ and $M$, 
then for any Bott tower \eqref{tower} with $M=B_n$, we have $2\mid h_j$ and $h_j^2=0$ for each $j$, and $\phi(c(P)) = c(M)$.
\end{lemma}

Note that for any holomorphic line bundle $Q$ on $B$, the projectivizations ${\mathbb P}(E)$ and ${\mathbb P}(E\otimes Q)$ are isomorphic to each other. In particular, for $B_j={\mathbb P}({\mathcal O}\oplus S_j)$ over $B_{j-1}$, one can replace $S_j$ by its dual $S_j^{\ast }$, as 
$$ {\mathcal O} \oplus S_j \cong (S_j^{\ast } \oplus {\mathcal O})\otimes S_j. $$
This replacement will not change $B_j$, but will affect the choice of  sections $L_j$ thus affecting $x_j$, while $h_j$ is replaced by $-h_j$. 

By the proof of the lemma above, we know that for any ${\mathbb Z}$-trivial Bott manifold $M^n$ and any Bott tower \eqref{tower} on 
$M$, if we write  $z_j=x_j+\frac{1}{2}h_j$, then $\{ z_1, \ldots , z_n\}$ is a set of   generators for the cohomology ring, with $z_j^2=0$ 
for each $j$. For any $2\leq j\leq n$, since $h_j^2=0$, we know that either $h_j=0$, or $h_j=  2q_j z_{\sigma (j)}$ for some positive 
integer $q_j$ and $\sigma (j) < j$. Here we used the fact that we can replace $S_j$ by $S_j^{\ast }$ without changing the Bott tower 
to ensure that these $q_j$ be positive.

From now on, we will make these choices, so $q_j>0$ whenever $h_j\neq 0$. That is, under our choices of these $S_j$, each 
$x_j$ is represented by the central sections of $\pi_j$, and each $z_j=x_j+ \frac{1}{2}h_j$ is represented by an effective divisor. We
have shown

\begin{lemma}
For any Bott tower \eqref{tower} we can choose the generator sets $\{ x_1, \ldots , x_n\}$ and $\{ z_1 , \ldots , z_n\}$ so that
(i) each $z_j^2=0$, (ii) each $z_j$ is represented by an effective divisor, and (iii) each $x_j$ is represented by a smooth hypersurface, 
which is iteself a Bott manifold of dimension $n-1$.
\end{lemma}

Obviously, for a given Bott manifold $M^n$, there are many Bott towers on it. So to sort out all distinct complex structures on $P=(S^2)^n$ 
given by the Bott manifolds, we need to find canonical representations for the ${\mathbb Z}$-trivial Bott manifolds. This is the goal of the 
next section.


\subsection{Bott diagrams.}

Let us denote by $A=\{ 1, 2, \ldots , n\}$ and write 
$A_0=\{ j\in A \mid h_j=0\}$. 
When $A_0\neq A$, we have a map $\sigma: A\setminus A_0 \rightarrow A$ satisfying $\sigma (j) < j$, given by the equation $h_j= 2q_j z_{\sigma (j)}$. Let us denote by $A_1=\sigma^{-1}(A_0)$, $A_2=\sigma^{-1}(A_1)$, and so on. It is easy to see that there exists some positive integer $r$ such that $A$ is the disjoint union of non-empty sets $A_0$, $A_1$, ... , $A_r$.

We will say that the {\em level} of $j\in A$ is $k$ if $j\in A_k$. It takes $\sigma$ exactly $k$ times to send a level $k$ element into $A_0$. 

\begin{definition}
For a given Bott tower \eqref{tower}, we define its {\em Bott diagram} to be the following data: each element of $A$ gives a vertex, each $j\in A \setminus A_0$ gives a vertical edge from $j$ to $\sigma (j)$, marked with a positive integer $q_j$. 
\end{definition}

In other words, a Bott diagram $G$ in dimension $n$ is a disjoint union $A=A_0\cup A_1\cup \cdots \cup A_r$ into $r+1$ nonempty subsets,  along with maps $A_r \rightarrow A_{r-1} \rightarrow \cdots \rightarrow  A_1 \rightarrow  A_0$ and a map $q:A\setminus A_0\rightarrow {\mathbb Z}^+$. Here $A$ is the set of $n$ elements and $r\geq 0$.

Two Bott diagrams are considered isomorphic, if there is a bijection from $A$ to $A$ which commutes with the partition of $A$ and the maps. 

We can arrange all the dots in $A_k$ at the same height, and will refer to that (imagined) horizontal line the level $k$ line (when $k=0$ we will 
also call it the base line). The diagram is a graph, with finitely many connected components which are {\it trees}. Each tree has a distinguished 
vertex, lying in $A_0$, which is the {\it root} of the tree. Thus from a combinatorial viewpoint, a Bott diagram is a {\it rooted forest}. Clearly, 
the Bott manifold is a product of lower dimensional ones, with each factor corresponds to a connected component of the Bott diagram. So 
$M^n$ is irreducible (in the sense that it is not the product of lower dimensional Bott manifolds) if and only if the Bott diagram is connected,
which occurs if and only if $A_0$ contains only one element. 

\begin{example*}
To illustrate how these diagrams work, let us first consider the case $n=2$. In this case we have only two possibilities for the Bott diagram: the first one just has two dots lying horizontally, with no edges, representing the surface ${\mathbb P}^1\times {\mathbb P}^1$; and the second one is two dots with a vertical edge connecting them, marked by a positive integer $q$. This corresponds to the Hirzebruch surface ${\mathbb F}_{2q}$.
\end{example*}

\begin{example*}
For $n=3$, we have three horizontal dots, corresponding to ${\mathbb P}^1\times {\mathbb P}^1 \times {\mathbb P}^1$; two dots on the base line, the third dot on top of the right one with a vertical edge marked mark $q$, corresponding to ${\mathbb P}^1 \times {\mathbb F}_{2q}$; one dot on the base line, two dots on the level $1$ line joining the base point by edges marked with $q$ and $p$, which corresponds to the fiber product $ {\mathbb F}_{2p}\times_{{\mathbb P}^1}  {\mathbb F}_{2q}$; and finally, we have three dots lined up in a vertical line, with two edges marked with $p$ and $q$ (with $p$ on top). In this case the threefold is $M^3={\mathbb P}( {\mathcal O}_B \oplus {\mathcal O}_B(-2p(C_0+qF)) ) $, where $B={\mathbb F}_{2q}$ is the Hirzebruch surface with $F$ the ruling and $C_0$ the central section (so $C_0^2=-2q$). 
\end{example*}


\subsection{Bott diagrams determine biholomorphism type.}  Our goal here is to complete the proof of our main theorem, by showing
the following:

\begin{theorem} \label{theorem}
Two Bott manifolds of dimension $n$ are biholomorphic to each other if and only if they have isomorphic Bott diagram. 
\end{theorem}
Since one can build up a Bott tower from the data of a Bott diagram, we just need to prove the ``only if" part of the statement, 
namely, if $f: M'\rightarrow M$ is a biholomorphism, then $M'$ and $M$ must have isomorphic Bott diagrams. 

Let us fix a Bott tower on $M$. By our previous discussion, we know that each $x_j$ is represented by a smooth hypersurface $X_j$ 
in the sense that $x_j=c_1(X_j)$, where the divisor $X_j$ is identified with the line bundle associated with it, and each $z_j$ is 
represented by an effective divisor. To be more precise, for any $j\in A_0$ of level $0$, $Z_j=X_j$ is irreducible. For any 
$j\in A_1$, $Z_j=X_j + q_jX_{\sigma (j)}$. For any $j\in A_2$, we have
$$ Z_j = X_j + q_j ( X_l + q_l X_{\sigma (l)}), \ \ \ \mbox{where} \ l=\sigma (j) \in A_1. $$
Note that each $X_j$ is itself a Bott manifold of dimension $n-1$, and the support of each $Z_j$ is a normal crossing divisor. 

In the case of a Bott tower, $X_1$ is a fiber of the composition map $\pi= \pi_2\circ \pi_3 \circ \cdots \circ \pi_n$ from $M^n$ to 
${\mathbb P}^1$, and $M$ is covered by the pencil of the divisor $X_1$. Given a Bott diagram, for any $j\in A_0$, we can 
choose a Bott tower on $M$ so that $j$ corresponds to the bottom layer, so we know that $M$ is covered by a pencil $|X_j|$ 
of the smooth hypersurface $X_j$, denote as $Y_t$ for $t\in {\mathbb P}^1$. These $Y_t$  do not intersect with each other. 

\vskip 5pt

We claim that, if $D$ is any effective divisor in $M$ homologous to $X_j$, then $D\in |X_j|$, namely, $D$ is a member of the pencil. 
To see this, first consider the special case when $D$ is irreducible. If $D$ is not in $|X_j|$, we may choose a member $Y$ in the 
pencil so that $D\cap Y \neq \phi$. This codimension $2$ subvariety is homologous to $0$, since $D\sim X_j$ and $x_j^2=0$, which 
is a contradiction since $M$ is projective. The same argument works when $D$ is just effective. 

\vskip 5pt

Now suppose $f: M' \rightarrow M$ is a biholomorphism. It induces a graded isomorphism $\phi = f^{\ast }$ between the cohomology rings. 
We want to show that $f$ induces an isomorphism between the Bott diagrams  as well. First we claim that $f$ induces a bijection 
between $A_0'$ and $A_0$, the set of level $0$ vertices. Let $\{ z_1, \ldots , z_n\}$ and $\{ z_1', \ldots , z_n'\}$ be generators on $M$ 
and $M'$ as before. Then $\phi (z_j)= z_{\tau (j)}'$ for some permutation $\tau$ on $A$ (the sign is positive since $z_j$, $z_j'$ are all 
represented by non-trivial effective divisors). If $j\in A_0$, then $M$ is covered by the pencil $X_j=Z_j$ of non-intersecting divisors. 
Consider the effective divisor $D=f(Z_{\tau (j)}')$ in $M$. $D$ represents $z_j$, thus is homologous to $Z_j$. By the claim above, 
we know that $D$ must be irreducible and is a member of the pencil $|X_j|$. This means that $\tau (j)$ lies in $A_0'$. So the level 
$0$ sets of $M$ and $M'$ are bijective to each other. 

Note that in the above argument, we furthermore obtained the fact that $f(X_{\tau (j)}') = X_j$ for $j\in A_0$. For $j\in A_0$, the smooth 
hypersurface $X_j$ is itself an $(n-1)$-dimensional Bott manifold. Its Bott diagram is obtained from that of $M^n$ by deleting the vertex 
corresponding to $j$, and pulling down one level in the tree above this vertex, while keeping everything else unchanged. We will call this 
new Bott diagram the {\em card} at vertex $j$. Now since $f(X_{\tau (j)}') = X_j$, the two Bott $(n-1)$-manifolds $X_{\tau (j)}'$ and $X_j$ 
are biholomorphic, so by induction on the dimension of the Bott manifolds, we see that the card of the Bott diagram $G'$ at vertex $\tau (j)$ 
must be isomorphic to the card of the Bott diagram $G$ at the vertex $j$. 

When $A_0$ has more than one elements, we have at least two cards, and we can use the set of cards to reconstruct the Bott diagram,
see Proposition \ref{reconstruction} and Remark \ref{reconstruction-labelled}. This implies that $G$ and $G'$ must be isomorphic to each other. When $A_0$ has only one element, the Bott diagrams agree as graphs, but we additionally need to show that the marking numbers $q_j$ for 
$j\in A_1$ should match those on $M'$ (see Remark \ref{reconstruction-labelled}).  

Without loss of generality, let us assume that $A_0=\{ 1\}=A_0'$. We already know that $f(X_1')=X_1$, and the card of the Bott graph $G$ at 
vertex $1$ is isomorphic to the card of the Bott graph $G'$ at vertex $1$. So $f$ gives bijection between $A_1$ and $A_1'$. Again without 
loss of generality, let us assume that $\phi (z_2) = z_2'$, where $2\in A_1$ and $2\in A_1'$. We have
$$ Z_2=X_2+ q_2X_1, \ \ Z_2' = X_2' + q_2'X_1'. $$
Consider the irreducible divisor $D=f(X_2')$. We have $D+q_2'X_1 \sim X_2 + q_2X_1$. If $D\neq X_2$, then $D\cap X_2$ is an effective 
cycle of codimension $2$ (could be trivial), and the intersection
$$ DX_2 + q_2' X_1X_2 \sim Z_2X_2 \sim Z_2(Z_2-q_2X_1) \sim -q_2 Z_2X_1 \sim -q_2 X_2X_1, $$
so the non-trivial effective cycle $D\cap X_2 + (q_2+q_2')X_1\cap X_2$ would be homologous to $0$ -- which is impossible since $M$ is 
projective. So $D$ must be equal to $X_2$, forcing $q_2'=q_2$. So the Bott diagrams of $M$ and $M'$ are indeed isomorphic to each other. 
This completes the proof of the theorem.

\section{Appendix: Reconstruction of rooted forests}\label{sec-reconstruction}

A famous problem in graph theory is the {\it reconstruction conjecture}. This conjecture asserts that finite graphs with 
at least three vertices are completely determined by their collection of vertex deleted graphs, see \cite{Harary}, \cite{BH}, \cite{Manvel}.
In this short appendix, we formulate and establish an analogue for rooted forests.

\begin{definition}
A contractable connected graph is called a {\it tree}. If such a graph is marked with a distinguished vertex (the {\it root}), 
we call it a {\it rooted tree}. If every connected component of a graph is a (rooted) tree, then we call the graph a (rooted)
forest. 
\end{definition}

Given a rooted tree $T$ with $n$ vertices, we can form an associated rooted forest with $n-1$ vertices by deleting the root $v$ of $T$
(and all incident edges). This leaves a forest with connected components $T_1, \ldots , T_k$, and we can pick a root on each tree $T_i$ 
to be the unique vertex $v_i$ of $T_i$ that was incident to $v$. We denote this rooted forest by $\hat T$, and call its individual trees
the {\it children} of the original tree $T$.

\begin{definition}
Let $F$ be a rooted forest, with connected components the rooted trees $T_1, \ldots , T_k$. Given a component $T_i$, we define
the associated {\it card} to be the rooted forest with components $T_1, \ldots, T_{i-1}, \hat T_i , T_{i+1}, \ldots ,T_k$, i.e. we replace 
the rooted tree $T_i$ by its collection of children. The forest $F$ has $k$ associated cards, each of which is a rooted forest.
\end{definition}

A graph is finite if it has finitely many vertices. We can now establish the reconstruction conjecture for rooted forests.

\begin{proposition}\label{reconstruction}
Let $F$ be a finite rooted forest. Then the set of cards of $F$ uniquely determines the forest $F$.
\end{proposition}

In other words, if one has a pair of forests $F_1, F_2$, and a bijection between the set of cards of $F_1$ and those of $F_2$, which
has the property that corresponding cards are isomorphic (as rooted forests), then the original forests have to be isomorphic.

\begin{proof}
We prove the statement using mathematical induction on the number of cards. 
Note that the number of cards coincides with the number of roots (and hence the number of connected components) in the original 
rooted forest $F$.

\vskip 5pt

\noindent \underline{Base case}: When there is only one card, we know that there is only one tree $T$ in the original forest $F$. The
individual trees in the single card are the children of $T$. We can thus reconstruct $T$ by taking a root vertex $v$, and for each of the 
children of $T$, connecting its root to the vertex $v$. The resulting rooted tree is the single tree in the forest $F$.

\vskip 5pt
 
\noindent \underline{Inductive step}: Let there be $n\geq 2$ cards in total. We run through the $n$ cards and locate a maximal tree (i.e. with
the maximal number of vertices) among all the trees appearing on all the cards. Of course there could be more than one such tree, but we 
pick one of them. Let us call this chosen maximal tree $T$. Our claim is that $T$ must be a rooted tree present in the original forest. If not, 
then $T$ appeared on the card after eliminating the root of one of the original trees $T_i$ of the forest. This means that $T$ is a proper 
subgraph of $T_i$, and that $T_i$ contains more vertices than $T$ (as the root $v_i$ of $T_i$ is not in $T$). Since $n\geq 2$, there is at least
one other card, arising from the deletion of another root $v_j$. The corresponding card contains $T_i$ as a rooted tree, contradicting 
the fact that $T$ was a maximal tree from all the cards. Note that this argument also shows that $T$ is not a child of any of the rooted trees 
in the original rooted forest.

Now that we have established $T$ is one of the original trees in the forest we are trying to rebuild, let us try and identify the multiplicity
with which it occurs in the forest. Assume the forest consists of $n$ rooted trees, and that $r$ of them are isomorphic to $T$ (where
$1\leq r \leq k$. Then there are precisely $r$ cards that contain $r-1$ copies of $T$, and $n-r$ cards that contain $r$ copies of $T$.
Thus, we may compute the integer $r$ from the set of cards.

Let $F'$ denote the forest obtained from the original forest $F$ by removing the $r$ copies of $T$. If we can reconstruct $F'$, then by
adding in $r$ copies of $T$, we will have reconstructed $F$. But note that the cards of $F'$ are easy to identify: just take the $n-r$
cards of $F$ that contain exactly $r$ copies of $T$, and remove from each of these cards the $r$ copies of $T$. The resulting $n-r$
rooted forests are the cards of $F'$. Since $r\geq 1$, the rooted forest $F'$ has $n-r<n$ cards, so by the inductive hypothesis, 
$F'$ can be reconstructed from its cards. Adding in $r$ disjoint copies of $T$ then produces $F$, and completes the proof of 
the Proposition.
\end{proof}

\begin{remark}\label{reconstruction-labelled}
Note that the proof of the proposition also holds for {\it labelled} rooted forests, where the cards are equipped with the natural 
induced labelling. In this setting, you need to additionally assume that the number of cards is $n\geq 2$ (i.e. this is the base case
of the induction, and is argued exactly like the inductive step above). When $n=1$, the only indeterminacy lies in the labels
for the edges in the rooted tree which are connected to the root vertex. These are obviously not recoverable from the single 
corresponding card. This is the reason for the additional argument at the end of the proof of Theorem \ref{theorem}.
\end{remark}



\end{document}